%

\def\astcirc{\circ}

%


\def\today{\ifcase\month\or January\or February\or
March\or April\or May\or June\or July\or August\or
September\or October\or November\or December\fi
\space\number\day, \number\year}




\def\dspace{\lineskip=2pt\baselineskip=18pt
\lineskiplimit=0pt}

\font \bbrm=cmbx10 at 12pt

\def\bigtype{\bbrm}

\hsize=13.5cm
\magnification=1200
\def\ce{\centerline}

\def\hb{\hfill\break}

\def\title #1{\null\bigskip\ce{\bigtype #1}
\bigskip}

\def\alp{\alpha}		
\def\bet{\beta}		
\def\gam{\gamma}		
\def\del{\delta}		
\def\eps{\varepsilon}

\def\kap{\kappa}
\def\lam{\lambda}		
\def\sig{\sigma}		

\def\ome{\omega}		


\def\calE{{\cal E}}

\def\calH{{\cal H}}

\def\calP{{\cal P}}



    
\font\tenboldgreek=cmmib10
 \font\sevenboldgreek=cmmib10 at 7pt
\font\fiveboldgreek=cmmib10 at 7pt
\newfam\bgfam
\textfont\bgfam=\tenboldgreek
\scriptfont\bgfam=\sevenboldgreek
\scriptscriptfont\bgfam=\fiveboldgreek

\mathchardef\ggarrow="7010

\font\tengerman=eufm10 \font\sevengerman=eufm7
\font\fivegerman=eufm5
\font\tendouble=msym10 \font\sevendouble=msym7
\font\fivedouble=msym5

\textfont4=\tengerman \scriptfont4=\sevengerman
\scriptscriptfont4=\fivegerman
\newfam\dbfam
\textfont\dbfam=\tendouble \scriptfont\dbfam=
\sevendouble
\scriptscriptfont\dbfam=\fivedouble

\mathchardef\ng="702D
\mathchardef\dbA="7041
\mathchardef\sm="7072
\mathchardef\nvdash="7030
\mathchardef\nldash="7031
\mathchardef\lne="7008
\mathchardef\sneq="7024
\mathchardef\spneq="7025
\mathchardef\sne="7028
\mathchardef\spne="7029
\mathchardef\ltms="706E
\mathchardef\tmsl="706F

\mathchardef\dbA="7041


\mathchardef\dbA="7041 
\mathchardef\dbB="7042 
\mathchardef\dbC="7043 
\mathchardef\dbD="7044 
\mathchardef\dbE="7045 
\mathchardef\dbF="7046 
\mathchardef\dbG="7047 
\mathchardef\dbH="7048 
\mathchardef\dbI="7049 
\mathchardef\dbJ="704A 
\mathchardef\dbK="704B 
\mathchardef\dbL="704C 
\mathchardef\dbM="704D 
\mathchardef\dbN="704E 
\mathchardef\dbO="704F 
\mathchardef\dbP="7050 
\mathchardef\dbQ="7051 
\mathchardef\dbR="7052 
\mathchardef\dbS="7053 
\mathchardef\dbT="7054 
\mathchardef\dbU="7055 
\mathchardef\dbV="7056 
\mathchardef\dbW="7057 
\mathchardef\dbX="7058 
\mathchardef\dbY="7059 
\mathchardef\dbZ="705A

\def\sdp{\times \hskip -0.3em {\raise 0.3ex
\hbox{$\scriptscriptstyle |$}}} 


\def\cf{{\rm \,cf\,}}

\def\dom{\mathop{\rm dom}\nolimits}

\def\min{\mathop{\rm min}}
\def\MOD{\mathop{\rm mod}}



\def\oB{{\overline B}}

\def\op{{\overline p}}


\def\omu{{\overline\mu}}



\def\udel{{\underline\del}}




\def\tilC{{\widetilde C}}

\def\tilf{{\widetilde f}}


\def\tilmu{{\widetilde\mu}}

\def\ddownarrow{\big\downarrow \hskip-0.70em\raise
2pt\hbox {$\big\downarrow$}}
\def\longright #1#2 {\smash{\mathop{\hbox to
#1pt {\rightarrowfill}}\limits_{#2}}}
\def\sqr#1#2{{\vcenter{\hrule height.#2pt\hbox{\vrule
width.#2pt height#1pt \kern#1pt \vrule width.#2pt}
\hrule height.#2pt}}}

\def\buildrul#1\under#2{\mathrel{\mathop{\null#2}
\limits_{#1}}}

\def\boxit#1{\vbox{\hrule\hbox{\vrule\kern3pt
\vbox{\kern3pt#1 \kern3pt}\kern3pt\vrule}\hrule}}

\def\suml{\sum\limits}
\def\prodl{\prod\limits}

\def\subheading#1{\medskip\goodbreak\noindent{\bf
#1.}\quad}

\def\sect#1{\goodbreak\bigskip\centerline{\bf#1}
\medskip}
\def\pr{\smallskip\noindent{\bf Proof:\quad}}
\def\onumber #1{\ooalign{\hfil\raise.07ex\hbox{
\hfill$\scriptstyle \,#1$\hfil}
\cr\cr{$\bigcirc$}}}
\def\onumber c{\ooalign{\hfil\raise.07ex\hbox
{\hfill$\scriptstyle \,c$\hfil}
\cr\cr{$\bigcirc$}}}
\def\alpcirc {\ooalign{\hfil\raise.07ex
\hbox{\hfill$\scriptstyle\alp\;$\hfill}\cr\cr
{$\bigcirc$}}}

\def\longmapright #1#2 {\smash{\mathop{\hbox to
#1pt {\rightarrowfill}}\limits^{#2}}}
\def\longmapleft #1 #2 {\smash{\mathop{\hbox to
#1 pt {\leftarrowfill}}\limits^{#2}}}

\def\references#1{\goodbreak\bigskip\par\centerline
{\bf References}\medskip\parindent=#1pt}
\def\ref#1{\par\smallskip\hang\indent\llap{\hbox
to \parindent{#1\hfil\enspace}}\ignorespaces}

\def\back{{\raise 2.5pt\hbox{$\,\scriptscriptstyle
\backslash\,$}}}
\def\bks{{\backslash}}
\def\part{\partial}
\def\lwr #1{\lower 5pt\hbox{$#1$}\hskip -3pt}
\def\rse #1{\hskip -3pt\raise 5pt\hbox{$#1$}}
\def\lwrs #1{\lower 4pt\hbox{$\scriptstyle #1$}
\hskip -2pt}
\def\rses #1{\hskip -2pt\raise 3pt\hbox
{$\scriptstyle #1$}}

\def\<#1{\left\langle{#1}\right\rangle}

\def\subinbn{{\subset\hskip-8pt\raise 0.95pt
\hbox{$\scriptscriptstyle\subset$}}}

\def\llvdash{\mathop{\|\hskip-2pt
\raise 3pt\hbox{\vrule height 0.25pt width 1.5cm}}}

\def\lvdash{\mathop{|\hskip-2pt \raise 3pt\hbox
{\vrule height 0.25pt width 1.5cm}}}

\def\fakebold#1{\leavevmode\setbox0=\hbox{#1}%
  \kern-.025em\copy0 \kern-\wd0
  \kern .025em\copy0 \kern-\wd0
  \kern-.025em\raise.0333em\box0 }

\font\msxmten=msxm10
\font\msxmseven=msxm7
\font\msxmfive=msxm5
\newfam\myfam
\textfont\myfam=\msxmten
\scriptfont\myfam=\msxmseven
\scriptscriptfont\myfam=\msxmfive
\mathchardef\rhookupone="7016
\mathchardef\ldh="700D
\mathchardef\leg="7053
\mathchardef\ANG="705E
\mathchardef\lcu="7070
\mathchardef\rcu="7071
\mathchardef\leseq="7035
\mathchardef\qeeg="703D
\mathchardef\qeel="7036
\mathchardef\blackbox="7004
\mathchardef\bbx="7003
\mathchardef\simsucc="7025

\def\rhookup{{\fam=\myfam \rhookupone}}

\def\bigsquare{{\fam=\myfam\bbx}}

\font\tencaps=cmcsc10
\def\smallcaps{\tencaps}

\def\author#1{\bigskip\ce{\smallcaps #1}\medskip}

\def\tagg{^{\prime\prime}}

\def\upddots{\mathinner{\mkern
1mu\raise 1pt \hbox{.}\mkern 2mu \mkern
2mu \raise 4pt\hbox{.}\mkern 1mu \raise 7pt\vbox
{\kern 7 pt\hbox{.}}} }

\def\varchi{\ooalign{{\raise
1.385pt\hbox{$\chi$}}\crcr\hbox{--}\crcr}}

\def\done{{1\hskip-2.5pt{\rm l}}}
\def\trianarrow{{\raise 2pt\hbox to 0.50cm
{\hrulefill}\triangleright}}

\input mssymb
\null
\overfullrule=0pt
\vskip2truecm
{\nopagenumbers
\overfullrule=0pt
\sect{MORE ON REAL-VALUED MEASURABLE CARDINALS}
\ce{\bf AND FORCING WITH IDEALS}
\bigskip
\ce{by}
$$\vbox{\halign{\tabskip2em\hfil#\hfil&\hfil#\hfil\cr
\bf Moti Gitik&\bf Saharon Shelah\cr
School of Mathematical Sciences&Hebrew University of
Jerusalem\cr
Sackler Faculty of Exact Sciences&Department of
Mathematics\cr
Tel Aviv University&Givat Ram, Jerusalem\cr
Ramat Aviv  69978 Israel&\cr}}$$

\dspace
{\narrower\medskip
\subheading{Abstract}
\item{(1)} It is shown that if $c$  is
real-valued measurable then the Maharam type of
$(c , \calP(c), \sig)$  is $2^c$.  This answers
a question of D.~Fremlin [Fr,(P2f)].
\smallskip
\item{(2)} A different construction of a model with a
real-valued measurable cardinal is given from
that of R.~Solovay [So]. This answers a question of
D.~Fremlin [Fr,(P1)]. 

\smallskip
\item{(3)} The forcing with a $\kap$-complete
ideal over a set $X$, $|X|\ge \kap$  cannot be
isomorphic to Random$\times$Cohen or Cohen$\times$Random.
 The result for $X=\kap$ was proved in [Gi-Sh1] but as
was pointed out to us by M.~Burke the application
of it in [Gi-Sh2] requires dealing with any $X$.

\vfill\eject}

\count0=1
\dspace
In Section 1 we deal with the Maharam types of
real-valued measurable cardinals.  The result (1) stated
in the abstract and its stronger version are
proved. The proofs are based on Shelah's strong
covering lemmas and his revised power set operation.

In Section 2 a model with a real-valued
measurable which is not obtained as the Solovay
one by forcing random reals over a model with a
measurable.

In Section 3, the result (3) stated in the abstract is
proved. 

Theorem 1.1 and the construction of Section 2 is
due to the first author.  
Theorem 1.2 is joint and the result of Section 3
is due to the second author.

We are grateful to David Fremlin for bringing
the questions on real-valued measurability to
our attention.  His excellent survey article
[Fr] gave the inspiration for the present
paper.  We wish to
thank the Max Burke for pointing out a missing stage in
the argument of [Gi-Sh2].\medskip}

\dspace
\def\llvdash{\mathop{\|\hskip-2pt \raise 3pt\hbox
{\vrule height 0.25pt width 0.50cm}}}

\sect{1.~~On Number of Cohen or Random Reals}

D.~Fremlin asked the following in [Fr,(P2f)]:

If $c$ is a real-valued measurable with
witnessing probability $\nu$, does it follow
that the Maharam type of $(c, \calP(c),
\nu)$  is $2^c$? or in equivalent formulation:\hb
If $c$ is a real-valued measurable does the
forcing with witnessing ideal isomorphic to the
forcing for adding $2^c$  random reals?   

The next theorem provides the affirmative
answer.

\proclaim Theorem 1.1.  Suppose that $I$  is a
$2^{\aleph_0}$-complete ideal over 
$2^{\aleph_0}$  and the forcing with it (i.e.
$\calP\big(2^{\aleph_0}\big)/I$) is isomorphic
to the adding of $\lam$-Cohen or $\lam$-random
reals.  Then $\lam=2^{2^{\aleph_0}}$.

\pr Suppose otherwise.  Denote $2^{\aleph_0}$
by $\kap$. Let $j:V\to N$  be a generic
elementary embedding.

\subheading{Claim 1}  $j(\kap)>(\lam^+)^V$.

\pr By a theorem of Prikry [Pr] (see also
[Gi-Sh2] for a generalization) for every $\tau
<\kap$  $2^\tau =2^{\aleph_0}=\kap$.  Then, in
$N$, $2^\kap =j(\kap)$.  But $(\calP(\kap))^V
\subseteq N$,  so $j(\kap)\ge (2^\kap)^V$.
By [Gi-Sh2], then $(2^\kap)^V=cov (\lam,
\kap, \aleph_1,2$.  So $cov(\lam, \kap,
\aleph_1, 2)\ge \lam^+$.  Clearly, 
$$cov (\lam,\kap, \aleph_1,2)\le cov (\lam,
\aleph_1,\aleph_1, 2)\le (cov (\lam,
\aleph_1,\aleph_1,2))^N\ .$$
The last inequality holds since $N$  is obtained
by a c.c.c. forcing and so every countable set
of ordinals in $N$  can be covered by a
countable set of $V$.  By Shelah [Sh430], in $N$
$cov (\lam, \aleph_1,\aleph_1,2)<j (\kap)$.
Hence $\lam^+\le cov (\lam,\kap,\aleph_1,2)\le
(cov (\lam,\aleph_1,\aleph, 2))^N<j(\kap)$.\hb
\hfill$\bigsquare$ of the claim. 

By Shelah [Sh 430] there is
$S\subseteq[\kap]^{\le\aleph_0}$ unbounded of
cardinality $\kap$.  For $\alp <\kap$  let
$S\rhookup\alp =\{ P\in S\big| P\subseteq\alp\}$.
Fix a function $f\in^\kap\kap$ representing $\kap$  is a
generic ultrapower and restrict everything to a condition
forcing this.

\subheading{Claim 2} 
$\{\alp <\kap \big|S\rhookup f(\alp)$  is
not unbounded in $[f(\alp)]^{\le\aleph_0}\}\in I$. 

\pr Otherwise, in a generic ultrapower $N$.
$j(S)\rhookup \kap =S$ will be bounded. I.e.
there will be some $t\subseteq \kap$  countable
such that for every $s\in S$ $s\not\supseteq t$.
Using c.c.c. of the forcing we find a countable subset
of $\kap$  in $V$, $t^*\supseteq t$. Since $S$
is unbounded in $V$ some $s\in S$ contains $t^*$. 
Contradiction.\hfill$\bigsquare$ of the claim.

Let $N$  be a generic ultrapower. By [Gi-Sh1]
there are in $N$  at least $\kap$  Cohen (or
random) reals over $V$.

\subheading{Claim 3}  There exists a sequence
$\langle r_\alp \mid \alp <\kap \rangle$  of
reals in $V$  so that 
\item{(1)} every real of $V$  appears in
$\langle r_\alp\mid\alp <\kap\rangle$.
\smallskip
\item{(2)} for almost all $\alp (\MOD I)$
$\langle r_{\alp +i}\mid i <f(\alp)\rangle$  are
$f(\alp)$-Cohen (random) generic over $L[ S\rhookup
f(\alp), \langle r_\bet\mid \bet <f(\alp)\rangle ]$.    

\pr Construct $\langle r_\alp \mid \alp
<\kap\rangle$  by induction.  On nonlimit stages
add reals in order to satisfy (1).  For limit
$\alp's$ with $S\rhookup f(\alp)$ unbounded in
$[f(\alp)]^{\le\aleph_0}$,  add $f(\alp)$-Cohen
(or random) reals. It is possible since there are at
least $\kap$ candidates in a generic ultrapower by
[Gi-Sh1].\hb
\hfill $\bigsquare$ of the claim.

Now work in $N$.  $rngf\rhookup A$ is unbounded
in $\kap$, for every $A\not\in I$. Let
$j(\langle r_\alp \mid\alp <\kap\rangle
)=\langle r_\alp \mid \alp <j(\kap)\rangle$
where $\langle r_\alp\mid\alp <\kap\rangle$  is
a sequence given by Claim 3.

Then, using Claim 3 in $N$ we can find some
$\alp^*<j(\kap)$  satisfying (2) of Claim 3 such
that $j(S)\rhookup\alp^*$  is unbounded in
$[\alp^*]^{\le \aleph_0}$  and $j(f)(\alp^*)\ge
(\lam^+)^V$.  It is possible since by Claim 1, $(\lam^+)^V
<j(\kap)$ and, in $V$  the range of $f$
restricted to a set not in $I$  is unbounded in
$\kap$. 

The following will provide the contradiction and
complete the proof of the theorem. 

\subheading{Claim 4} $\langle r_{\alp^*+i}\mid
i<j(f)(\alp^*)\rangle$  is a sequence of Cohen
(random) reals over $V$.

\pr $\langle r_{\alp^*+i}\mid
i<j(f)(\alp^*)\rangle$ is Cohen (random)-generic
over $L[j(S)\rhookup j(f)(\alp^*),\ \langle
r_\bet\mid\bet <j(f)(\alp^*)\rangle ]$.
In particular, also over
$L[j(S)\rhookup j(f)(\alp^*),\ \langle r_\bet \mid\bet
<\kap \rangle ]$.  But $\langle r_\bet \mid \bet
<\kap\rangle$  is a list of all the reals of $V$.
It is enough to show that $j(S)\rhookup j(f)(\alp^*)\cap
([j(f)(\alp^*)]^{\le\aleph_0})^V$ is unbounded.

We note that the very strong covering holds between
$L[S\langle r_\alp\mid \alp <\kap\rangle]$  and
$L[j(S), \langle r_\alp \mid \alp <\kap\rangle ]$.
Also, for any $b\in V\cap H_\kap$,  the same is
true about $L[b,S,\langle r_\alp \mid \alp
<\kap]\rangle ]$ and $L[b,j(S),\langle r_\alp
\mid \alp <\kap \rangle ]$.  Then the same holds
in $N$, between
$$L[b,j(S)\rhookup j(f)(\alp^*),\ \langle r_\bet
\mid \bet <j(f)(\alp^*)\rangle ]\ {\rm and}\ L[b,j(S),\
\langle r_\bet \mid\bet <j(\kap )\rangle ]\ ,$$
for every $b\in V\cap H_\kap$.  Now consider $V$,
$V[G]$. Using c.c.c. it is easy to find the
winning strategy for $V$.  Take a play which is
according such a strategy and of the length $<\kap$,
in sense of [Sh420, 2.6, $\mu (3)]$.  Then it
belongs to $L[j(S),\langle r_\bet \mid \bet
<j(\kap )\rangle ]$  since the last model agrees
with $V[G]$  about small sequences.  Now use the
Shelah Strong Covering [Sh580] for $L[j(S)\rhookup
j(f)(\alp^*),\langle r_\bet\mid\bet
<j(f)(\alp^*)]$, $V$  and $V[G]$.  We will obtain 
$$P\in V\cap L[j(S)\rhookup j(f)(\alp^*),\ \langle
r_\bet \mid \bet <j(f)(\alp^*)]\cap [j(f)(\alp^*)]^{\le
\aleph_0}\ .$$
Since $|P|=\aleph_0$  all its subsets in $V$  are
also in $L[S,\langle r_\alp |\alp <\kap\rangle ]\subseteq
L[j(S)\rhookup j(f)(\alp^*)$,  $\langle r_\bet \mid
\bet <j(f)(\alp^*)\rangle ]$.  Hence the Cohen
(random) genericity over the last model is
equivalent the Cohen (random) genericity over
$V$.\hfill$\bigsquare$

Let us now prove a stronger statement which
relies on a different property.

\proclaim Theorem 1.2.  Suppose that $I$  is a
$\kap$-complete ideal over $\kap$  and the
forcing with it (i.e. $\calP(\kap)/I$) is
isomorphic to the forcing for adding
$\lam$-Cohen or $\lam$-random reals.  Assume also\hb 
(*) some condition forces that ``$j(\kap)\ge
(2^\kap)^V$", where $j$  is a generic embedding.
Then $\lam >2^{<\kap}$  implies $\lam =2^\kap$.

\pr Without loss of generality, let us assume that the
weakest condition, i.e. $\kap$  forces (*).  Suppose
that $\lam <2^\kap$.  Then $\lam < j(\kap)$.  By
[Sh430], then in a generic ultrapower $N$ $cov
(\lam , \aleph_1, \aleph_1,2)<j(\kap)$.  However, in
$V$  $\lam^+\le cov (\lam,\kap,\aleph_1,2)\le
cov (\lam,\aleph_1\aleph_1,2))^N<j(\kap)$.
The first inequality holds by [Gi-Sh2].  Hence
$\lam^+<j(\kap)$.

Now, by [Sh 430, 460, 2.6], there are regular
$\del <\mu <\kap$ such that $cov (\lam, \mu, \mu,\del)=
\lam$.

Let us assume for simplification of the notation
that $\mu =\aleph_2$, $\del =\aleph_1$.

\subheading{Claim 1}  There is a sequence of
reals $\langle r_\alp \mid\alp <\lam^+\rangle$
in a generic ultrapower such that for every
$s\subseteq \ome_1$  the final segment of
$\langle r_\alp \mid \alp <\lam^+\rangle$  are
Cohen (or random) generic over $L[s]$. 

\pr Let $N$  be a generic ultrapower.  Then
${}^\kap\!N\subseteq N$  where ${}^\kap\!N$ is in the
sense of the generic extension.  First note that if
$\langle s_\alp \mid\alp<\kap\rangle$ is a sequence of
$\kap$ Cohen (random) reals over $V$,  then it belongs
to $N$.  Clearly every $s\subseteq\ome_1$  in $N$  (or
the same in $V[G]$)  is a name in $V$
interpreted using $\aleph_1$  Cohen (random) reals only.
Hence the final segment of $\langle s_\alp\mid\alp
<\kap\rangle$  will be generic over $L[s]$.
Then, in $V$, for every regular $\del <\kap$
there will be a sequence $\langle t_\alp
\mid\alp <\del\rangle$ such that for every $s\subseteq
\ome_1$  its final segment is Cohen (random) generic
over $L[s]$.  Back in $N$,  we use this for $\del =\lam^+$
which is still below $j(\kap)$.\hb
\hfill$\bigsquare$ of the claim.

Let us fix such a sequence $\langle r_\alp\mid\alp
<\lam^+\rangle$ in $N$.  We split it into blocks
each of the length $\ome_1$.  Denote such
changed sequence by $\langle r_{\alp i}\mid\alp <\lam^+,
i<\ome_1\rangle$.  Now back in $V$,
let us use the fact that $cov(\lam,\aleph_2,\aleph_2,
\aleph_1)=\lam$.  We know that for every $\alp
<\lam^+$  the block
$\langle{\vtop{\offinterlineskip\hbox
{$r$}\hbox{$\scriptscriptstyle\sim\alp_i$}}}\mid i<\ome_1
\rangle$ is added by using only $\ome_1$  Cohen (or
random) reals from the $\lam$  Cohen (or Random)
reals, $\langle C_\bet\mid\bet <\lam\rangle$.
More precisely, for every $\alp <\lam^+$  there
is $t_\alp :{}^{\ome_>}\ome_1\longrightarrow {}^{\ome_>}
\ome_1$ in $V$  and a set of indexes $b_\alp\in
V[G]$  and some enumeration of it $\langle
\xi_{\alp i}\mid i<\ome_1\rangle$ such that
$\langle r_{\alp i}\mid i<\ome_1\rangle$  is the
image of $\langle C_{\xi_{\alp_i}}\mid
i<\ome_1\rangle$ under $t_\alp$.  Since $2^{<\kap}<\lam$,
we can assume w.l. of g.  that for some $t$  for
every $\alp <\lam^+$  $t_\alp =t$.  We view
$b_\alp$  as $\langle \langle \xi_{\alp i},i\rangle\mid
i<\ome_1\rangle$ i.e. a subset of $\lam\times \ome_1$  of
cardinality $\aleph_1$.  Since $|\lam|=|\lam \times \ome_1|$
we still can use $cov (\lam, \ome_2,\ome_2,\ome_1)=\lam$
and find $b\subseteq \lam\times\ome_1$,
$|b|=\aleph_1$ such for $\lam^+$  $\alp's$
$|b\cap b_\alp |=\aleph_1$  and $t$ applied
to $b$  provides a nontrivial information about
$\langle r_{\alp i}\mid i<\ome_1\rangle$.  But
now the final segment of $\langle r_{\alp
i}\mid\alp <\lam^+,\ i<\ome_1\rangle$  cannot be
Cohen (random) generic over $L[t,b]$.
Contradiction.\hfill$\bigsquare$

\sect{2.~~Another Construction of a Model with a
Real-Valued Measurable Cardinal}

In this section we construct a model with a
real-valued measurable cardinal which differs
from the Solovay original.  This answers
negatively a question of D.~Fremlin\hb
[Fr, (P1)]:

Let $N$  be a model of ZFC and $\kap \in N$  a
real-valued measurable cardinal in $N$.  Does it
follow that there are inner models $M\subseteq N$  such
that $\kap$  is a measurable in $M$  and
$M$-generic filter $G$  for a random real p.o.
set over $M$  such that $G\in N$  and $N\cap
\calP(\kap)\subseteq M[G]$.

Suppose that $\kap$ is a measurable and $GCH$
holds.  We define a forcing notion $P$  as
follows:

\subheading{Definition 2.1}  $P$ consists of all
triples $p=\langle p_0,p_1,p_2\rangle$  so that
\item{(1)} $p_0\subseteq\kap$
\smallskip
\item{(2)} $p_1$  is a function with domain
contained in $p_0$
\smallskip
\item{(3)} $p_2$  is a function defined over
inaccessibles $\le\kap$
\smallskip
\item{(4)} for every inaccessible $\del$ $|p_0\cap\del
|<\del$,  $|\dom p_1\cap \del|<\del$  and $|\dom
p_2\cap\del |<\del$
\smallskip
\item{(5)} for every $\alp\in\dom p_1$
$p_1(\alp)\subseteq \alp$
\smallskip
\item{(6)} every element of $p_0$  is an ordinal
of cofinality $\aleph_0$
\smallskip
\item{(7)} for every limit ordinal $\bet$ if
$\cf\bet >\aleph_0$,  then $p_0\cap\bet$  is not
stationary in $\bet$  and if $\cf\bet=\aleph_0$
then $\bet\bks (p\cap \bet)$  is unbounded in
$\bet$
\smallskip
\item{(8)} for every $\alp\in\dom p_2$
$p_2(\alp)$  is a closed subset of $\alp$
disjoint with $p_0$.

\subheading{Definition 2.2} Let $p,q\in P$
$p=\langle p_0,p_1,p_2\rangle$  and $q=\langle
q_0,q_1,q_2\rangle$.  Then $p\ge q$  iff  
\item{(1)} $p_1\subseteq q_1$
\smallskip
\item{(2)} $\dom p_2\supseteq\dom q_2$  and for
every $\alp\in\dom q_2$  $p_2(\alp)$  is an end
extension of $q_2(\alp)$
\smallskip
\item{(3)} $p_0\supseteq q_0$
\smallskip
\item{(4)} for every $\del <\kap$,  $\del$ is an
inaccessible or a limit of inaccessibles and
$\del^*$  is the least inaccessible above $\del$
then $p_0\cap [\del,\del^*)$  is an end
extension of $q_0\cap[\del,\del^*)$. 

The forcing $P$  is intended to add three
objects.  Thus, the first coordinates of $P$ are
producing a subset $S$  of $\kap$ which is
stationary in $V[S]$  and reflecting only in
inaccessibles.  The second coordinate is
responsible for a kind of diamond sequence  
over $S$  and the last coordinate adds clubs
preventing reflection of $S$  at inaccessibles
and its stationarity. 

The forcing $P$  destroys the measurability of
$\kap$  once used over $V=L[\mu]$.  It is bad
for our purpose.  We are going to use a certain
subforcing of $P$  which will preserve
measurability and contain the projection of $P$
to the first two coordinates.  But first let us
study basic properties of $P$.

Let $P_0=\{ p_0\mid\exists \langle p_1,p_2\rangle$
$\langle p_0,p_1,p_2\rangle\in P\},\ P_{01}=\langle
\langle p_0,p_1\rangle |\exists p_2\ \langle
p_0,p_1,p_2\rangle\in P\}$.  Let $\alp$  be an
inaccessible.  We denote by $P\rhookup\alp$  the
set
$$\{ \langle p_0\cap\alp,\ p_1\rhookup\alp,\
p_2\rhookup\alp\rangle\mid \langle
p_0,p_1,p_2\rangle\in P\}$$ 
and by $P\bks\alp$  the set
$$\{\langle p_0\bks\alp,\ p_2\rhookup
[\alp,\kap),\ p_2\rhookup [\alp +1,\kap)\rangle
\mid \langle p_0,p_1,p_2\rangle \in P\}\ ,$$
$P_0\rhookup\alp, P_{01}\rhookup \alp$ and $P_0\bks\alp$,
$P_{01}\bks\alp$  are defined similarly.

The following is standard.

\subheading{Claim 2.3}  Let $\alp$  be an
inaccessible then the following holds
\item{(1)} $P =P\rhookup\alp\times P\bks\alp$
\smallskip
\item{(2)} $P_0=P_0\rhookup\alp\times
P_0\bks\alp$
\smallskip
\item{(3)} $P_{01}=P_{01}\rhookup\alp\times
P_{01}\bks\alp$.
\item{}

\smallskip
Let $\alp <\kap$  be a limit ordinal and $Q$ a
forcing notion.

Consider the following game Game $(Q,\alp)$:
\input pictex
$$\matrix{{\bf I}&q_1&&&&q_3&\cdots&&&\cdots&&\cr
&&{\beginpicture
\setcoordinatesystem units <0.25cm,0.25cm>
\plot  9.335  5.683  8.668  5.905 /
\plot  8.668  5.905  8.890  5.207 /
\plot  8.446  5.683  8.668  4.985 /
\linethickness=0pt
\putrectangle corners at  8.446  5.905 and  9.335  4.985
\endpicture}
&&{\beginpicture
\setcoordinatesystem units < 0.25cm, 0.25cm>
\plot  8.287  6.636  8.064  5.969 /
\plot  8.064  5.969  8.763  6.160 /
\plot  8.287  5.747  8.985  5.937 /
\linethickness=0pt
\putrectangle corners at  8.064  6.636 and  8.985  5.747
\endpicture}&&&&\ge&&&\cr
{\bf II}&&&q_2&&&&\cdots&&q_\bet&\cdots&q_\alp\cr}$$
where Players I, II are building an increasing
sequence of elements of $Q$, I at even stages
and II at odds.  If at some stage $\bet <\alp$
II cannot continue i.e. there is no $q$  above
$\{ q_\bet'|\bet'<\bet\}$  then I wins.
Otherwise II wins. 

\subheading{Claim 2.4}  The player II has a winning
strategy in the game Game $(P\bks\alp,\alp^+)$
for every inaccessible $\alp$.

\pr Let $\alp$  be an inaccessible.  We define a
winning strategy $\sig$  for Player II in the
Game $(P\bks\alp,\alp^+)$.

Let $\del >\alp$  be an inaccessible but not limit one. 
Denote by $\del^-$ the supremum of inaccessibles below
$\del$. 

Let $p\in P\bks\alp$. We define $\op$ to be the
condition obtained from $p=\langle p_0, p_1,p_2\rangle$
by adding $\sup (p_0\cap[\del^-,\del)) +\sup(p_2(\del))$
to $p_2(\del)$ if $p_0\cap [\del^-,\del)
\not=\emptyset$ or $p_0\cap [\del^-,\del)=\emptyset$ but
$p_0\cap\del^-$  is unbounded in $\del^-$,
where $\del$ runs over inaccessibles above
$\alp$  which are not limit inaccessibles and
$\sup (p_2(\del))=0$  whenever $\del\not\in\dom p_2$.

Now we define $\sig$  to be dependent only on
the last move of $I$ at successive stages of the game.
Set $\sig (p_{\bet +1})=\op_{\bet+1}$.  If $\bet
\le\alp^+$ is limit and the game up to $\bet$
$$\matrix{p_1&&p_3&\cdots&&&\cdots\cr
&p_2&&\cdots&&p_\gam&\cr}$$
was played according to $\sig$.
Then set $\sig (\langle p_\gam\mid\gam
<\bet\rangle)=$ the closure of $\bigcup_{\gam <\bet}
p_\gam$.  More precisely, let $\sig (\langle p_\gam
\mid\gam <\bet \rangle)=\langle p^0,p^1,p^2\rangle$ 
where $p^0=\bigcup_{\gam <\bet}p^0_\gam$,
$p^1=\bigcup_{\gam <\bet}p^1_\gam$  and
$\dom (p^2)=\bigcup_{\gam <\bet}\dom (p^2_\gam)$, 
$p^2(\xi)=\bigcup\{p^2_\gam(\xi)\mid\gam <\bet,\xi\in\dom
p^2_\gam\}\cup\{\sup (\cup\{p^2_\gam (\xi)\mid\xi\in\dom
p^\gam_2,\gam <\bet\})\}$,for $\xi\in\dom p^2$.

We need to check that such defined $p$  is a
condition.  The only problem is to show that
$p^0$  does not reflect at any $\tau$,
$\aleph_0<cf\tau <\tau$.  So let $\tau$  be an
ordinal such that $\aleph_0<cf\tau <\tau$  and
$p^0\cap\tau$  is unbounded in $\tau$.  Pick
$\del$  to be the first inaccessible above
$\tau$.  Then $\del^-\le\tau$.  If $\del^-<\tau$,
then starting with some $\gam_0<\bet\ $
$p^0_\gam\cap[\del^-,\del)\not=\emptyset$.  But
then $p^2(\del)$  will be a club of $\tau$ disjoint
to $p^0\cap\tau$.  Suppose now that
$\del^-=\tau$.  Then $C_0=\{\sup\big(\bigcup_{\gam'<\gam}
(p^0_{\gam'}\cap\tau)\big)\mid\gam <\bet ,\ \gam\ {\rm
limit}\}$ is a club of $\tau$.  Also
$C_1=\{\xi <\tau\mid\xi\ \hbox{is a limit of
inaccessibles}\}$ is a club of $\tau$. Let
$\xi\in C_0\cap C_1$.  Then for some
inaccessible $\udel <\tau$  $\xi =\udel^-$. Let
$\gam_\xi <\bet$  be so that
$\xi=\sup\bigcup_{\gam'<\gam_\xi}(p^0_{\gam'}\cap\tau)$.
Then, $\xi\in p^2_{\gam_\xi}(\udel)$,  by the
definition of $\sig$ and the procedure. Thus,
this provides a club of $\tau$ disjoint to
$p^0_\bet\cap\tau$.\hfill $\bigsquare$ of the
claim.

The following is now trivial.

\proclaim Claim 2.5.  $P$  preserves
cofinalities and does not add new functions from
ordinals less than the first inaccessible into
$V$.

Let $U$ be a normal measure over $\kap$  and $j:V\to N$ 
the corresponding elementary embedding.  Then,
in $N$, $j(P)=j(P)\rhookup\kap\times j(P)\bks
\kap$.  Clearly, $(j(P)\rhookup \kap)^N=P$.  Now
let us produce inside $V$ an $N$-generic subset
of $j(P)\bks\kap$  with the set over the first
coordinate nonstationary in $V$.  

\subheading{Claim 2.6} There exists $\langle S,
F_1, F_2\rangle$  such that 
\item{(a)} $\langle S, F_1,F_2\rangle$  is
$j(P)\bks\kap^+$ generic over $N$ 
\item{(b)} $S$ is not stationary subset of
$j(\kap)$
\item{(c)} $S$  does not reflect.
\item{}

\pr Let $\langle D_\alp |\alp <\kap^+\rangle$  be
the list of dense open subsets of $j(P)\bks\kap$
of $N$.  Let $\sig\in N$ be a winning strategy
for Player II in Game $(j(P)\bks\kap, \kap^+)$.
It exists by Claim 2.4 applied in $N$ to $j(P)$.
Play the game from $V$  so that I plays at stage
$\bet +1$  an element $P_{\bet +1}$ of $D_\bet$ which
is above $p_\bet$,  where $\bet<\kap^+$.  We will finish
with a desired $N$-generic set.\hfill$\bigsquare$  

Force with $P$ over $V$. Let $G$  be a generic
subset.  We denote $\bigcup \{ p_0\mid\exists
\langle p_1,p_2\rangle\ \langle p_0,p_1,p_2\rangle$
$\in G\}$
by $S$. For every $\alp\in S$  let $A_\alp
=\bigcup\{ p_1(\alp)\mid\exists\langle
p_0,p_1,p_2\rangle\in G$ and $\alp\in\dom p_1\}$
and for inaccessible $\del\le\kap$  let $C_\del
=\bigcup \{ p_2(\del)\mid \exists \langle
p_0,p_1,p_2\rangle\in G$ and $\del\in\dom p_2\}$.
Then $S\subseteq\kap$, and for every inaccessible
$\del\le\kap$  $C_\del$  is a club of $\del$
disjoint to $S$.   

\subheading{Claim 2.7}  $S$  is a stationary
nonreflecting subset of $\kap$  in $V[S,\langle
A_\alp\mid\alp\in S\rangle$, $\langle C_\del\mid
\del\ {\rm inaccessible\ and}\ \del<\kap\rangle]$.

\pr $\langle C_\del\mid\del <\kap\rangle$  are
witnessing the nonreflection.  Suppose that $S$
is nonstationary in
$V[S,\langle A_\alp\mid\alp\in S\rangle$, 
$\langle C_\del\mid\del\ {\rm inaccessible\
and}\ \del <\kap\rangle ]$.  Return back to $V$
and work with names. Suppose for simplicity that
the empty condition forces the nonstationarity
of $S$.  Let $\buildrul\sim\under C$  be a name of
witnessing club.  Pick an elementary submodel $N$ of
$V_{2^{2^\kap}}$ such that $P,\buildrul\sim\under C\in
N,\ |N|<\kap$  and $N\cap\kap$  is an ordinal of
cofinality $\aleph_0$.  Let $\alp =N\cap\kap$  and
$\langle \alp_n\mid n<\ome\rangle$  be a cofinal in $\alp$
sequence.  Now by induction we construct an increasing
sequence $\langle p_i\mid i<\ome\rangle$  of
conditions of $P\rhookup \kap$  (i.e. $P$
without the information on a club of $\kap$
disjoint to $S$) such that for every $i<\ome$  
\item{(a)} $p_i\in N$
\item{(b)} $p_i$  decides the first element of
$\buildrul\sim\under C$  above $\alp_i$
\item{(c)} $\sup (p_i)_0\ge\alp_i$.
\item{}
\noindent
where $p_i=\langle (p_i)_0,(p_i)_1, (p_i)_2\rangle$
\smallskip
Now, in $V$,  let
$$p=\langle\bigcup_{i<\ome} (p_i)_0\cup
\{\alp\},\bigcup_{i<\ome}(p_i)_1,\
\{\langle\del,\cup \{ (p_i)_2(\del)\mid
i<\ome,\ \del\in\dom (p_i)_2\}\rangle\ .$$
Then $p\llvdash\alp\in\buildrul\sim\under C
\cap\buildrul\sim\under S$.
Contradiction.\hfill$\bigsquare$ 

\subheading{Claim 2.8}  $\kap$  is a measurable
cardinal in $V[S,\langle A_\alp\mid\alp\in S\rangle ]$.

\pr Just note that in $N$  $j(P_{01})=P_{01}\times
j(P_{01})\bks\kap^+$, since nothing is done in
the interval $[\kap,\kap^+]$  by this forcing.
By Claim 2.6, there is a $j(P_{01})\bks\kap^+$
generic over $N$ set in $V$.  Thus it is easy to
extend $j$  to the embedding of $V[S,\langle
A_\alp\mid\alp\in S\rangle ]$.  This insures
the measurability of $\kap$.\hfill$\bigsquare$

Notice, that $\kap$  is not measurable in
$V[S,\langle A_\alp\mid \alp\in S\rangle ,\ \langle
C_\del\mid\del <\kap ,\ \del\ {\rm inaccessible}\rangle]$ 
since $S$ is a stationary nonreflecting subset of $\kap$.

Now, over a model $V[S,\langle A_\alp\mid\alp\in
S\rangle ]$ we are going to force a Boolean
algebra $B$  such that:
\item{(a)} $\kap$  is still measurable in $V[S,\langle
A_\alp\mid\in S\rangle , B]$  
\smallskip
\item{(b)} $j^*(B)\big/ G(B)$  is isomorphic to the
adding of $j(\kap)$-Random reals, where $j^*$
the elementary embedding of $V[S,\langle
A_\alp\mid\alp\in S\rangle, B]$ into its
ultrapower and $G(B)$  a generic subset of $B$.
\item{}

\smallskip
First let us review some basics of product
measure algebras.  We refer to D.~Fremlin [Fr2]
for detailed presentation.

Suppose that $B$ is a {$\sig$-algebra\/}, i.e. a
Boolean algebra all of whose countable suprema
exist.  A {\it measure on\/} $B$  is a function
$\mu :B\to [0,1]$  so that:   (a) $\mu (1_B)=1$,
and (b) whenever $\{ b_n\mid n\in\ome\}\subseteq
B$  with $b_n\wedge b_n=0$ for $n\not= m$,  then
$\mu\big(V_nb_n\big)=\suml_n\mu (b_n)$.
If in addition $\mu$  is {\it positive\/} (i.e.
$\mu(b)=0$  {\it iff\/} $b=0$),  then we say
that $\langle B,\mu\rangle$ is a {\it measure
algebra\/}.  A measure algebra is always a
complete Boolean algebra.

Suppose now that $I$ is a set, and $\langle
B_i,\mu_i\rangle$  for $i\in I$ are measure
algebras.  Call $C\in\prodl_{i\eps I}B_i$ a
{\it cylinder iff\/} $C(i)$  is the unit element
of $B_i$, except for a finite number of
coordinates $i$. Let $B\supseteq\prodl_{i\eps
I}B_i$  be the $\sig$-algebra generated by the
cylinders.  It is known that there is a unique
measure $\mu$  on $B$ so that $\mu(C)=\prodl_{i\eps I}
\mu_i(C(i))$ for any cylinder $C$.  $\mu$  may not be 
positive, but there is a standard strategy: Let
$I=\{ b\in B\mid \mu (b)=0\}$.  Then $I$ is an
ideal, and $\oB=B/I$  as usual is a
$\sig$-algebra consisting of equivalence classes
$[b]$  for $b\in B$  (where $[b]=[c]$ {\it iff\/} the
symmetric difference $(b-c)\vee (c-b)\eps I)$.
We can define a positive measure $\omu$  on
$\oB$  by: $\omu([b])=\mu(b)$.  Thus, $\langle
\oB,\omu\rangle$ is a measure algebra, called
the product {\it measure algebra\/} of the $\langle
B_i,\mu_i\rangle$'s. 

Let ${\bf 2}$  be the basic measure algebra $\langle
P(2),\mu\rangle$  where $\mu$  is the measure:
$\mu(\emptyset)= 0,\ \mu (\{ 0\})=\mu(\{
1\})={1\over 2}$,  and $\mu(\{ 0,1\})=1$.  For
any set $I$, let ${\bf 2}^I$  denote the product
measure algebra of $I$  copies of ${\bf 2}$.  We can
then force with ${\bf 2}^I$  with the natural
proviso:  $b$  is a {\it stronger\/} condition
than $c$  {\it iff\/} $0<b\le c$ in ${\bf 2}^I$.
This forcing obviously has the $\ome_1$-c.c. 

For $I=\ome$, ${\bf 2}^I$  is just the usual random
real forcing and for $I=\lam$ $\ {\bf 2}^I$  is the
$\lam$-random real forcing.  Let us denote them
by Random and Random($\lam$) respectively.

We consider the $\sig$-algebras $B_\alp\subseteq
\prodl_{i<\alp}(\calP(2))_i$ generated by the
cylinders, where $\alp\le\kap$  and $(\calP(2))_i$
is just $i$-th copy of $\calP(2)$.  The desired
algebra $B$ will be $B_\kap/I_\kap$,  where the
ideal $I_\kap$ of ``null" sets is going to be
added generically.  More precisely for $\alp$'s
of countable cofinality $I_\alp$'s will be added
by forcing and then for $\bet\le\kap$  of
uncountable cofinality $I_\bet$  will the union
of $I_\alp$'s, where $\alp <\bet$,  $cf\alp =\aleph_0$.  
The sequence of ideals $\langle
I_\alp\mid\alp\le\kap\rangle$  will be in $V_1=V[S,\langle
A_\alp\mid\alp\in S\rangle$,  $\langle
C_\del\mid\del <\kap ,\ \del\ {\rm inaccessible}\rangle]$.

Let us work in $V_1$.  We define by induction on
$\alp <\kap$  a measure $\mu_\alp$ on $B_\alp$.
Then $I_\alp$  will be the ideal of $\mu_\alp$-measure
null sets.  Actually there will be a lot of different
measures over $B_\alp$'s.  We would like to
prevent $B_\kap$ (and even its subalgebras of
power $\kap$) from carrying a measure.
For this purpose, the ``diamond" sequence $\langle
A_\alp\mid\alp\in S\rangle$  will be used to destroy
possible candidates. 

If $\alp <\min S$,  then let $\mu_\alp$  be the
usual product measure over $B_\alp$, i.e. one
generated by attaching weight $1/2$  to $\{0\}$
and $\{ 1\}$,  $0$ to $\emptyset$ and $1$ to $\{
0,1\}$ in every component $(\calP(2))_i$
$(i<\alp)$  of the product $\prodl_{i<\alp}(\calP(2))_i$. 
Set $I_\alp=\{ X\in B_\alp\mid\mu_\alp(X)=0\}$.

Suppose now that $\alp <\kap$  and for every
$\bet <\alp$  the measure $\mu_\bet$  over $B_\bet$
was already defined.  We need to define $\mu_\alp$
over $B_\alp$.

\subheading{Case 1} $\alp\notin S$.\hb 
Pick an increasing continuous sequence $\langle
\alp_\tau\mid\tau <cf\alp\rangle$  witnessing
nonstationarity of $S\cap\alp$.  In case,
$cf\alp =\aleph_0$  just use $\ome$-sequence 
unbounded in $\alp$  and disjoint with $S$.
For every $\tau <cf\alp$ let $\mu(\tau)$  be
$\mu_{\alp_{\tau +1}}\rhookup\big(B_{\alp_{\tau+1}}
\rhookup [\alp_\tau,\alp_{\tau +1})\big)$, where
$B_{\alp_{\tau +1}}\rhookup
[\alp_\tau,\alp_{\tau+1})$  is the subalgebra of
$\prodl_{\alp_\tau \le i<\alp_{\tau +1}}(\calP(2)))_i$ 
generated by the cylinders.

Let $\mu_\alp$  be the product measure of $\langle
< B_{\alp_{\tau +1}}\rhookup [\alp_\tau,\alp_{\tau
+1})$, $\mu(\tau)>\mid \tau <cf \alp\rangle$.

Notice that $\alp_\tau\notin S$.  Therefore by
induction, we can assume that for a limit $\tau$
the measure $\mu_{\alp_\tau}$  over
$B_{\alp_\tau}$  is the product measure of
$\langle < B_{\alp_\nu}\rhookup [\alp_\nu,
\alp_{\nu+1})$, $\mu (\nu)>\mid\nu <\tau\rangle$. 

\subheading{Case 2} $\alp\in S$.\hb
Suppose that $A_\alp$  codes in some reasonable
fashion sequences $\langle \alp_n\mid <\ome\rangle$,
$\langle\varphi_n\mid n<\ome \rangle$,
$\langle\mu (n)\mid n<\ome\rangle$  and $\langle
a_n\mid n<\ome\rangle$  so that for every $n
<\ome$  
\item{(a)} $\langle\alp_n\mid n<\ome\rangle$  is
a cofinal in $\alp$ sequence
\smallskip
\item{(b)} $a_n$  is a countable subset of
$[\alp_n,\alp_{n+1})$ 
\smallskip
\item{(c)} $\mu(n)$  is a measure over
$B_{\alp_{n+1}}\rhookup a_n$  respecting the
ideal $I_{\alp_{n+1}}\rhookup a_n$, i.e. for
every $X\in B_{\alp_{n+1}}\rhookup a_n$  $\mu
(n)(X)=0$ iff $X\in I_{\alp_{n+1}}$
\smallskip
\item{(d)} $\varphi_n$: Random $\leftrightarrow
B_{\alp_{n+1}}\rhookup a_n$  is a measure
algebra isomorphism.

Denote $B_{\alp_{n+1}}\rhookup a_n$  by $B(n)$.
Let us define a measure $\tilmu (n)$  over $B(n)$.
Thus for every $n<\ome$  let us change the value
$\mu (n)(\varphi_n(\{ 0\} ))$  from $1/2$  to
$1-{1\over\pi^2 n^2}$ and those of $\varphi_n(\{
1\} )$  from $1/2$ to ${1\over\pi^2n^2}$.  Let
$\tilmu (n)$  be the measure obtained from $\mu (n)$ in
such a fashion.  Clearly, such local changes
have no effect on the set of measure zero.  Namely,
for every $X\in B(n)$  $\mu(n)(X)=0$  iff
$\tilmu (n)(X)=0$.

Define now the measure $\mu_\alp$  over $B_\alp$
as the product measure of the measure algebras
$\langle B(n),\tilmu (n)\rangle$ $(n<\ome)$
together with all the rest, i.e.
$$\langle B_{\alp_{n+1}}\rhookup
([\alp_n,\alp_{n+1})\bks a_n),\ \mu_{\alp_{n+1}}
\rhookup ([\alp_n,\alp_{n+1})\bks a_n)\rangle\ .$$

We claim that $\varphi =\bigcup_{n<\ome}\varphi_n$ 
cannot be extended to complete embedding into
$\langle B_\alp,\mu_\alp\rangle$.  The reason is
that under $\varphi$  the measure of the set
$\bigcap_{n<\ome}\varphi_n(\{ 0\} )$  should be
zero, but $\mu_\alp\big(\bigcap_{n<\ome}\varphi_n
(\{0\})\big)=\prodl_{n<\ome}\tilmu(n)\big(\varphi_n(\{
0\})\big)=\prodl_{n<\ome}\big(1-{1\over\pi^2n^2}\big)$
which equals sin$(1)\not= 0$  by the Euler formula. 

Notice, that the ideal $I_\alp$  of
$\mu_\alp$-measure zero sets will not be
effected if for finitely many $n$'s the measures
$\mu(n)$  will be used in the product instead of
$\tilmu(n)$'s.  Also, if in the previous
construction we will do everything above some
$\alp_{n_o}$  for fixed $n_0<\ome$, i.e. 
we will define the measure over $B_\alp\rhookup
[\alp_{n_0},\alp)$  instead of all $B_\alp$ call
it $\mu_\alp\rhookup [\alp_{n_0},\alp)$  and its
ideal $I_\alp\rhookup [\alp_{n_0},\alp)$,  then,
for every $X\in B_\alp$  $X\in I_\alp$ iff
$X\rhookup \alp_{n_0}\in I_{\alp_{n_0}}$ and
$X\rhookup [\alp_{n_0},\alp)\in I_\alp \rhookup
[\alp_{n_0},\alp)$.  This means that once having
$I_{\alp_n}$'s, initial segments of measures
$\langle\tilmu (n)\mid n<\ome\rangle$ have no
effect on $I_\alp$.  
This observation will be crucial further for
showing measurability of $\kap$.

If $A_\alp$  does not guess the sequences as
above, then we proceed as in Case 1. 

This completes the definition of $\langle\mu_\alp\mid\alp
<\kap\rangle$  and hence also $\langle I_\alp\mid\alp\le
\kap\rangle$.

We set $B=B_\kap\big/I_\kap$.
Let $V_2=V[S,\langle A_\alp\mid \alp\in
S\rangle$,  $\langle I_\alp\mid\alp
<\kap,\cf\alp =\aleph_0\rangle]$.  Then for
every $\alp<\kap$  $I_\alp\in V_2$.  So $B\in
V_2$. We will show the following claim which has
a proof similar to 2.7.

\subheading{Claim 2.9}  Random$(\kap)$  does not
embed into $B$  in $V_1$  and also in $V_2$.  

\pr Notice that $V_2$ and $V_1$ have the same
reals.  So if $\varphi$  is an embedding of
Random$(\kap)$ into $B$ in $V_2$  then $\varphi$
will be also such embedding in $V_1$.  Hence let
us prove the claim for $V_1$.

Suppose otherwise.  Let $\varphi:
{\rm Random}(\kap)\longrightarrow B$  witnessing
embedding.  Back in $V$  let us work with names.
Let $\buildrul\sim\under\varphi$ be a name of $\varphi$
and assume for simplicity that the empty condition forces
this.

Pick $N$ and $\langle\alp_n\mid n<\ome\rangle$
to be as in Claim 2.7 with $\buildrul\sim\under\varphi$
replacing $\tilC$.  

We define sequences of conditions of $P\rhookup
\kap$ of $\{ p_n\mid n<\ome\}\subseteq N$,  of
ordinals $\langle\bet_n\mid n<\ome\rangle$,
countable sets $\langle a_n\mid n<\ome\rangle$
and embedding $\langle\varphi_n\mid n<\ome\rangle$ so that
\item{(a)} $\sup (p_n)_0\ge\alp_n$
\item{(b)} $\bet_n\ge\alp_n$
\item{(c)} $p_n\llvdash\tagg\buildrul\sim\under\varphi(\{
0\}\bet_n)$ has nontrivial intersection with
$\buildrul\sim\under B\rhookup [\bet_n,\bet_{n+1})\tagg$,
where $\{ 0\}_{\bet_n}\in (\calP(2))_{\bet_n}$ i.e.
the $\bet_n$-th copy of $\calP(2)$.
\item{(d)} $a_n\subseteq\bet_n$
\item{(e)} $\varphi_n$ embeds $(\calP(2))_{\bet_n}$
into $B_\kap\rhookup a_n$
\item{(f)} $p_n\llvdash\tagg\check\varphi_n$ is
equal to $\buildrul\sim\under\varphi\rhookup
(\calP(2))_{\bet_n}\MOD{\vtop{\offinterlineskip\hbox
{$I_\kap$}\hbox{$\scriptscriptstyle\sim$}}}\tagg$

Since the forcing does not add new countable
sequences of elements of $V$, there is no
problem in carrying out the induction.

Denote by $\mu(n)$  the measure over
$B_\kap\rhookup a_n$  induced by $\varphi_n$.

Now let $A_\alp\subseteq \alp$  be a code for
such sequences 
$$\langle \bet_n\mid n<\ome\rangle,\ \langle
a_n\mid n<\ome\rangle,\ \langle\varphi_n\mid
n<\ome\rangle$$
and $\langle\mu (n)\mid n<\ome\rangle$.

Set $p=\langle\bigcup_{n<\ome}(p_n)_0\cup\{\alp\},
\bigcup_{n<\ome}(p_n)_1\cup\{ \langle \alp, A_\alp
\rangle \}$,
$$\{\langle\del , \cup\{(p_n)_2(\del)\mid
n<\ome,\del\in\dom (p_n)_2\}\rangle$$ 
Then $p\llvdash\tagg\buildrul\sim\under \varphi$ does
not embed
$\prodl_{n<\ome}(\calP(2))_{\bet_n}$  into
$B_\alp\big/ {\vtop{\offinterlineskip\hbox
{$I_\alp$}\hbox{$\scriptscriptstyle\sim$}}}\tagg$,
by the definition of $I_\alp\tagg$,.  Hence, also
$p\llvdash\tagg\buildrul\sim\under\varphi$ does not embed
Random$(\kap)$ into $B\tagg$. Contradiction.\hfill
$\bigsquare$

\subheading{Claim 2.10}  $\kap$  is a measurable
cardinal in $V_1$.

\pr Let $j:V\longrightarrow N$  be an elementary
embedding witnessing the measurability of
$\kap$.  We like to extend it to an embedding 
$$j^*:V[S,\langle A_\alp\mid\alp\in S\rangle ,\
\langle I_\alp\mid\alp <\kap\ {\rm and}\ \cf
\alp =\aleph_0\rangle]$$
$$\longrightarrow N[S^*,\ \langle
A_\alp\mid\alp\in S^*\rangle ,\ \langle
I_\alp\mid\alp <j(\kap)\ {\rm and}\ \cf\alp
=\aleph_0\rangle]\ .$$

By Claim 2.6, $j$  extends to
$$j': V[S,\langle A_\alp\mid \alp\in S\rangle ]
\longrightarrow N[S^*,\langle A_\alp\mid\alp\in
S^*\rangle ]$$ 
where $\langle S^*\bks S,\ \langle A_\alp\mid\alp\in
S^*\bks S >\rangle$ $\in V$  is $j(P_{01})\bks\kap$
generic over $N$.  We like to produce ideals
$\langle I_\alp\rhookup (B_{j(\kap)}\rhookup
[\kap, j(\kap))\mid \kap <\alp <j (\kap),
\cf\alp =\aleph_0\rangle$ generic over $N$  but
in $V$. In order to define $\langle
I_\alp\mid\alp <\kap,\cf \alp =\aleph_0\rangle$
we used clubs witnessing nonreflection of $S$,
i.e. $\langle C_\del\mid\del <\kap ,\ \del\
{\rm inaccessible}\rangle$.  By Claim 2.6, the
only club which is needed in order to extend $j$
but is missing in $V$ is $C_\kap$.  But, we
define generically $I_\alp$'s only for $\alp$'s
of cofinality $\aleph_0$  and moreover initial
segments have no influence on such $I_\alp$'s.
This means that the definition of $\langle
I_\alp\rhookup B_{j(\kap)}\rhookup [\kap,j(\kap))\mid\kap
<\alp <j(\kap), cf\alp=\aleph_0\rangle$ can be
carried out completely inside $N[S^*\bks
S,\langle A_\alp\mid\alp\in S^*\bks S\rangle$,
$\langle C_\del\mid\kap <\del <j(\kap), \del\
{\rm inaccessible}\rangle ]$. All the sets
$S^*\bks S,\langle A_\alp\mid\alp\in S^*\bks
S\rangle$ and $\langle C_\del\mid\kap <\del
<j(\kap),\ \del\ {\rm inaccessible\ of}\
N\rangle$ can be found inside $V$  by Claim
2.6.  Hence we have enough sets to extend $j$
to $j^*$.  Thus, the measurability of $\kap$  is
preserved in $V_2=V[S,\langle A_\alp\mid\alp\in
S\rangle$, $\langle I_\alp\mid\alp <\kap$  and
$\cf\alp =\aleph_0\rangle ]$.
\hfill$\bigsquare$

Let $j^*:V_2\longrightarrow N_2=N[S^*,\langle
A_\alp\mid\alp\in S^*\rangle\ ,\langle I\alp\mid
\alp <j(\kap)$ and $\cf\alp =\aleph_0\rangle ]$
be the embedding of Claim 2.10.

\subheading{Claim 2.11} $j^*(B)\rhookup [\kap,
j(\kap))$  is isomorphic to Random$(\kap^+)$ in
$V_2$.

\pr By Claim 2.6, there are in $V$  and hence in
$V_2$ clubs $\langle C_\del\mid \kap <\del\le
j(\kap),\ \del\ \hbox{is\ an}$
${\rm inaccessible\ in}$ $N\rangle$
witnessing nonreflection of $S^*\bks S$  in
every $N$-inaccessible $\del\le j(\kap)$.  Using
them we define measures $\mu_\alp$  over $B_\alp
\rhookup [\kap, \alp)$ agreeing with ideals
$I_\alp$  for every $\alp$, $\kap <\alp \le
j(\kap)$  as it was done for $B_\alp$'s below
$\kap$  in $V_1$.  The final measure $\mu_{j(\kap)}$
will turn $B_{j(\kap)}\rhookup [\kap, j(\kap))$
into measure algebra.  Since $|j(\kap)|=\kap^+$,
by Maharam theorem, see [Fr2] we obtain the
desired result.\hfill$\bigsquare$ 

So we have the following:

\proclaim Theorem 2.  $\kap$  be the measurable
cardinal of $L[\mu]=V$  (the minimal model with
a measurable).  Then $\kap$  is a real valued
measurable in $V^B_2$  but for every submodel
$V'$  of $V_2^B$  if $\kap$  is a measurable in
$V'$,  then there is no $G\in V^B_2$  which is
Random$(\kap)$  generic over $V'$.

\pr Suppose that $V'\subseteq V_2^B$,  $\kap$  is
a measurable in $V'$, $G\subseteq$ Random$(\kap)$
generic over $V'$ and $G\in V^B_2$.  But then
$G$  is also Random$(\kap)$  generic over $V=L[\mu]$,
since $L[\mu]\subseteq V'$  by its minimality.
But $V$ and $V_2$  have the some countable
sequences of ordinals.  So, $G$  will be
Random$(\kap)$-generic also over $V_2$.  This
means that Random$(\kap)$  embeds $B$,  which is
impossible by Claim 2.9.\hfill$\bigsquare$
 
\sect{3.~~The Forcing with Ideal Cannot be
Isomorphic to Cohen$\fakebold{$\times$}$Random}
\vskip-0.25truecm
\ce{\bf or Random$\fakebold{$\times$}$Cohen}

The result for $\kap$-complete ideals over
$\kap$  was proved in [Gi-Sh1].  Max Bruke
pointed out that the application of this in
[Gi-Sh2] requires the result also for less than
$\kap$ complete ideals as well.  The purpose of
this section is to close this gap. 

\proclaim Theorem 3.1.  Suppose that $I$  is a
$\ome_1$-complete ideal over some $\kap$ then
the forcing with ideal (i.e. $\calP(\kap)/I$)
cannot be isomorphic to Cohen$\times$Random or
Random$\times$Cohen. 

\pr Let us deal with Random$\times$Cohen case.
The Cohen$\times$Random case is similar. 

Suppose otherwise.  $\calP(\kap)/I\simeq$
Random$\times$Cohen.  Without loss of generality
for some $\kap_1\le \kap$  and $f:\kap\to\kap_1$
$\kap\llvdash_{\calP(\kap)/I}{\rm "}\kap_1$  is the
critical point of the generic embedding and $\tilf$ 
represents $\kap_1$  in the ultrapower". 
Define an ideal $J$  over $\kap_1$  to be the
set of all $A\subseteq\kap_1$  such that
$f^{-1\prime\prime} (A)\in I$.  Denote
$Q=\calP(\kap)/I$  and $Q_1=\calP(\kap_1)/J$.
Then $Q_1$  is a complete subordering of $Q$.
  We define a $Q_1$-name 
$\buildrul\sim\under\tau=\{\eta\in {}^{\ome >}\!2\mid$
the condition $\left(\done_{\rm Random},\eta\right)$ is
compatible with every element of $\buildrul\sim\under
G(Q_1)\}$.  For $\eta,\nu\in{}^{\ome >}\!2$ let
us write $\eta\triangleright\nu$  if the
sequence $\eta$  extends the sequence $\nu$. 
The following two claims are obvious.

\subheading{Claim 3.2} $\buildrul\sim\under\tau$  is a
$Q_1$-name of a nonempty subset of ${}^{\ome >}\!2$ 
closed under initial segments with no
$\triangleleft$-maximal element and hence a tree.

\subheading{Claim 3.3} $\llvdash_Q\tagg$ the Cohen 
real is an $\ome$-branch of $\buildrul\sim\under\tau$".

\subheading{Claim 3.4}  There is no $p\in Q_1$
and $\eta\in {}^{\ome>}\!2$  such that $p\llvdash_{Q_1}$
``for every $v\in {}^{\ome >}\!2\ v\triangleright\eta$ 
implies $v\in\buildrul\sim\under\tau$".

\pr Suppose otherwise.
Let $p,\eta$ be witnessing this. Then above $p$
the forcing notion $Q_1$  is a complete
subordering of Random.  But it has to add a
real.  Hence it is isomorphic to Random which  
is impossible by [Gi-Sh1].\hfill$\bigsquare$

Let $T^*=\{ T\mid T\subseteq {}^{\ome>}\!2$  is
a tree and for every $n<\ome$, $\eta\in {}^n\!2$
there are $v\triangleright\eta$ and $k < \ome$  such that $v\in
{}^{k}\!2$  and $v\notin T\}$.
Consider also $T^*_m=\{T\cap {}^m\!2\mid T\in
T^*\}$  for $m<\ome$.  $T^*$  can be viewed as a
tree if we identify it with
$\bigcup\limits_{m<\ome}T^*_m$ and define an order by
setting $t_1\triangleleft t_2$  iff for some
$m<\ome$  $t_1=t_2\cap{}^{m>}\!2$.  Then, clearly,
$$\llvdash_{Q_1}\tagg\buildrul\sim\under\tau\ 
\hbox{is an}\ \ome -{\rm branch\ of}\ T^*{}\!\tagg\ .$$  

\subheading{Claim 3.5} Suppose that $n<\ome$,
$q_0\in$ Random, $\eta\in {}^n\!2$.  Then there
are $m<\ome$, $q, v_0, v_1,t_0, t_1$  such that 
\item{(a)} $q\in$ Random and $q\ge q_0$.
\smallskip
\item{(b)} $\eta \triangleleft v_0,v_1\in
{}^{\ome >}\!2$
\smallskip
\item{(c)} $t_0, t_1\subseteq {}^{m\ge}\! 2$
and $t_0\not= t_1$.
\smallskip
\item{(d)} $(q,v_i)\llvdash\buildrul\sim\under\tau
\cap {}^{m\ge}\!2 =t_i$ for $i<2$.

\pr Find first some $q'\ge q_0$ and
$v_0\triangleleft\eta$  deciding $\buildrul\sim\under\tau
\cap {}^{m\ge}\!2$.  Let $t_0$  be the decided value,
i.e. $(q',v_0)\llvdash\buildrul\sim\under\tau\cap{}^{m\ge}
2=t_0$.  By the Claim 3.4 there will be $m<\ome$
and $v\triangleleft\eta$, $v\in {}^m\!2\bks t_0$.
Find some $(q,v_1)\ge (q',v)$  deciding $\buildrul\sim
\under\tau\cap {}^m\!2$.  Let $t_1$  be the
forced value, i.e. $(q,v_1)\llvdash\buildrul\sim\under
\tau\cap{}^m\!2=t_1$.  Since $(q,v_1)\llvdash v_1\in    
\buildrul\sim\under\tau$,  we have $(q,v_1)\llvdash
v=v_1\rhookup m\in\buildrul\sim\under\tau$.  But
this means $t_o\not= t_1$.\hfill$\bigsquare$ 

\subheading{Claim 3.6}  Suppose that $n,k<\ome$
and $q_0\in$ Random.  Then there are $q\in$
Random, $m<\ome$, $\langle v_{\eta,\ell}\mid \eta\in
{}^n\! 2,\ \ell <k\rangle$ and $\langle
t_{\eta,\ell}\mid\eta\in {}^n\! 2$,  $\ell
<k\rangle$ such that
\item{(a)} $q\ge q_0$
\smallskip
\item{(b)} $m\ge n$
\smallskip
\item{(c)} for every $\eta_1,\eta_2\in {}^n\!
2,\ \ell_1,\ell_2 <k$ $t_{\eta_1,\ell_1}=t_{\eta_2,
\ell_2}$ iff $(\eta_1,\ell_1)=(\eta_2,\ell_2)$
\smallskip
\item{(d)} for every $\eta\in {}^n\!2,\ \ell <k$
$$\eta\triangleleft v_{\eta,\ell}\in {}^{\ome
>}\!2\ ,\ t_{\eta,\ell}\in T^*_m\ {\rm and}\ 
(q,v_{\eta,\ell})\llvdash \buildrul\sim\under\tau\cap
{}^m\!2=t_{\eta,\ell}\ .$$

\pr Just use the previous claim enough times.
Thus, first, we generate a tree of $k\cdot (2^n+1)$ 
possibilities for one $\eta\in {}^n\! 2$
and then we repeat the argument of Claim 3.5 on all
$\eta$'s.\hfill$\bigsquare$

\subheading{Claim 3.7}  For every $n<\ome$, $k<\ome,
q'\in$ Random and $\calE >0$  there are $m<\ome$,
$q\ge q', \{ q_\ell \mid \ell < \ell^*\}\subseteq$ Random
pairwise disjoint, $\langle v_{\eta,
\ell,j}\mid\eta\in {}^n\! 2,\ell<\ell^*\ ,\ j<k\rangle$
and $\langle t_{\eta,\ell,j}\mid\eta\in {}^n\!2,\
\ell <\ell^*,\ j<k\rangle$ such that  
\smallskip
\item{(a)} $Lb(q)\ge 1-\calE$ ($Lb$ denotes the
Lebesgue measure)
\smallskip
\item{(b)} $q=\bigcup\limits_{\ell <\ell^*}q_\ell$
\smallskip
\item{(c)} if $\eta\in{}^n\!2$,  $\ell <\ell^*$
and $j<k$  then $v_{\eta,\ell, j}\in {}^{\ome >}\! 2$ 
$v_{\eta,\ell,j}\triangleright\eta$  and
$(q_\ell,v_{\eta,\ell,j})\llvdash\buildrul\sim\under\tau
\cap{}^{m\ge}\! 2=t_{\eta,\ell,j}$.
\smallskip
\item{(d)} for every $\ell <\ell^*$
$$t_{\eta_1\ell, j_1}=t_{\eta_2,\ell, j_2}\quad
{\rm iff}\quad (\eta_1,j_2)=(\eta_2,j_2)\ .$$  

\pr We define by induction $q_\ell$'s using
Claim 3.6.  Thus if $\langle q_i\mid
i\le\ell\rangle$ is defined then we apply Claim
3.6 to ${}^\ome\!2\bks\bigcup\limits_{i\le \ell}q_i$.
The process stops after we reach $\ell^*$ s.t.
$Lb(\bigcup\limits_{\ell <\ell^*}q_\ell)\ge
1-\calE$.\hfill$\bigsquare$

\subheading{Claim 3.8} For every $n<\ome$  and
$\calE >0$  there are $m, n\le m<\ome$  and a function
$H:T^*_m\longrightarrow 2$ such that for every
$\eta\in {}^n\!2$ and $i\in 2$  we can find
$q^{i,\eta}, \ell^{i,\eta} <\ome,\langle q_\ell^{i,\eta}
\mid\ell <\ell^{i,\eta}\rangle$ and $\langle
v_\ell^{i,\eta}\mid\ell <\ell^{i,\eta}\rangle$
such that for every $i<2$  and $\ell <\ell^{i,\eta}$ 
\item{(a)} $\eta\triangleleft v_\ell^{i,\eta}\in
{}^{\ome >}\!2$
\smallskip
\item{(b)} $q^{i,\eta},q_\ell^{i,\eta}\in$
Random and $q^{i,\eta}=\bigcup\limits_{\ell
<\ell^i}q_\ell^{i,\eta}$
\smallskip
\item{(c)} $Lb(q^{i,\eta})\ge 1-\calE$
\smallskip
\item{(d)} $\langle q_\ell^{i,\eta}\mid\ell
<\ell^{i,\eta}\rangle$ are pairwise disjoint
\smallskip
\item{(e)} $\left(q^{i,\eta}_\ell,v_\ell^{i,\eta}\right)
\llvdash$ ``$H\left(
\buildrul\sim\under\tau\cap{}^{m\ge}\!2\right)=i$".

\pr For every $\eta\in {}^n\!2$  and $t\in
T_m^*$ let $I_{\eta, t}=\left\{q\in {\rm
Random}\mid\ \hbox{there is}\ v\in{}^{\ome
>}\!2,\ v\triangleleft\eta\right.$\hb
such that $\left.(q,v)\llvdash\tagg
\buildrul\sim\under\tau \cap{}^{m\ge}\!2=t\tagg
\right\}$.  Let $\left\{ q_{\eta,t,\ell}\mid\ell
<\ell_{\eta,t}\le \ome\right\}$ be a maximal antichain
subset of $I_{\eta,t}$.  Let $q^*_{\eta,t}=
\bigcup\limits_{\ell <\ell_{\eta,t}}q_{\eta,t,\ell}$. 
Then $\bigcup\limits_{t\in T^*_m}q^*_{\eta,t}={}^\ome\!2$
mod null set, since $\left\{ q_{\eta,t,\ell}\mid
t\in T^*_m\right.$, $\left.\ell <\ell_{\eta,t}\right\}$ is a
predense subset of Random (but not necessarily antichain).
So $Lb\big(\bigcup_{t\in T^*_m}q^*_{\eta, t}\big)=1$.  

It is enough to prove the following statement:\hb
(*) There exists $H:T_m^*\longrightarrow 2$ so
that for every $\eta\in {}^n\!2$  and $i<2$
$$Lb\left(\bigcup\{ q^*_{\eta,t}\mid t\in T^*_m\
{\rm and}\ H(t)=i\}\right)\ge 1-{\calE\over 2}\ .$$
Since then we will be able to find a maximal
antichain $\langle q^i_\ell\mid\ell<\ell^*\le\ome\rangle$
in Random above $\bigcup\left\{q^*_{\eta,t}\mid t\in
T^*_m\ {\rm and}\ H(t)=i\right\}$ together with $\langle
v_\ell^i\mid \ell <\ell^*\rangle$  and $\langle
t^i_\ell\mid \ell <\ell^*\rangle$ so that
$$\left(q^i_\ell,v^i_\ell\right)\llvdash\tagg
\buildrul\sim\under\tau\cap{}^{m\ge}\!2=t^i_\ell\
{\rm and}\ H(t_\ell)=i\tagg\ .$$
In order to reduce $\ell^*$ to a finite $\ell^i$
we note that the precision here is $1-{\calE\over 2}$
but only $1-\calE$ is needed. 

So let us prove (*).  We consider the set $\calH$  of all
functions $H:T_m^*\longrightarrow 2$.  It is
finite but more transparent is to look at it as a
probability space.  All $H\in\calH$  with the
same probability.  So we choose $H(t)\in \{
0,1\}$  independently for the $t\in T^*_m$  with
probability $1/2$. 

We use $m$ given by Claim 3.7 for our $n,\calE'$
much smaller than $\calE$  and $k$  large
enough.  Given $\eta\in {}^n\! 2$ and $i<2$.  We
consider the probability of 
$$\left(Lb\Big(\bigcup\left\{q^*_{\eta,t}\mid
t\in T^*_m\ ,\ H(t)=i\right\}\Big)\ge
1-{\calE\over 2}\right)$$
in $\calH$.  It is $\le 1$, as the value is
always $\le 1$ and is $\ge 1-{1\over 2^k}$.  In
order to prove the last inequality, let us use
$\{ q_\ell\mid\ell <\ell^*\}$ of Claim 3.7.
Thus 
$$\eqalign{&Lb\Big(\bigcup \Big\{ q^*_{\eta,t}\mid t\in
T^*_m\ {\rm and}\ H(t)=i\Big\}\cr  
&=\suml_{\ell <\ell^*}Lb \Big(q_\ell \cap\bigcup\Big\{
q^*_{\eta,t}\mid t\in T^*_m\ {\rm and}\
H(t)=i\Big\}\Big)+Lb\Big[\Big({}^\ome\!2\bks
\bigcup\limits_{\ell<\ell^*}q_\ell\Big)\cap\cr
&\Big(\bigcup\Big\{q^*_{\eta,t}\mid t\in
T^*_m\ {\rm and}\ H(t)=i\Big\}\Big)\Big]\ge\cr
&\ge\suml_{\ell <\ell^*}\Big[Lb(q_\ell)\times
\Big(Lb\Big(q_\ell\cap\bigcup\Big\{q^*_{\eta,t}\mid
t\in T^*_m\ {\rm and}\ H(t)=i\Big\}\Big)\Big]
\big/Lb(q_\ell)-\calE'\cr}$$
Since (a) of 3.7, $Lb\left({}^\ome\!2\bks
\bigcup\limits_{\ell<\ell^*}q_\ell\right)<\calE'$.  Now
it suffices to show that for each $\ell
<\ell^*$ 
$$Lb\left(q_\ell\cap\cup\left\{q^*_{\eta,t}\mid
t\in T^*_m\ ,\ H(t)=i\right\}\right)\big/Lb(q_\ell)\ge
1-{\calE\over 4}$$
holds for enough $H$'s.
But $v_{\eta,\ell,j}$, $t_{\eta,\ell,j}$ $(j<k)$
of 3.7 are witnessing that the probability in
$\calH$  of the failure is $\le {1\over 2^k}$.
Just in order to fail, $H$  should take the
value $1-i$ on $t_{\eta, \ell,j}$  for every
$j<k$  and the probability of $0,1$ are equal.
The probability of the failure for some $\eta\in
{}^n\!2, i\in 2$  is then $\le{2^{n+1}\over 2^k}$.
So, picking $k$  large enough comparatively to
$n$  we will insure that most $H\in\calH$  are
fine, whereas we need only one.\hfill$\bigsquare$  

Now using Claim 3.8, we define by induction on
$j<\ome$  $\ n_j,m_j, H_j, \langle q^{i,\eta}_j\mid
i<2, \eta\in^{n_j}2\rangle$ $\langle
q^{i,\eta}_{j,\ell}\mid\ell <\ell_j$ $i<2,
\eta\in{}^{nj}\!2\rangle$ and $\langle
v_{j,\ell}^{i,\eta}\mid\ell <\ell_j,
i<2,\eta\in {}^{n_j}\!2\rangle$  such that  
\item{(1)} $n_j <m_j<n_{j+1}$
\smallskip
\item{(2)} $m_j H_j$, $\langle q_j^{i,\eta}\mid i<2,
\eta\in {}^{n_j}\!2\rangle, \langle q^{i,\eta}_{j,\ell}
\mid\ell <\ell_j^{\eta, i}, i<2, \eta\in {}^{n_j}\!2
\rangle\ {\rm and}\ \langle v_{j,\ell}^{i,\eta}\mid \ell
<\ell_j^{i,\eta}, i <2, \eta\in{}^{n_j}\!2\rangle$   
are given by Claim 3.8 for $n=n_j$  and $\calE={1\over
2^{2^{n_j}}}$
\smallskip
\item{(3)} length $\left(v_{j,\ell}^{i,\eta}\right)
<n_{j+1}$ for every $i<2$,  $\eta\in{}^{n_j}\!2,\ell
<\ell_j^{i,\eta}$. 

Now define a $Q_1$-name $\buildrul\sim\under\sig\in
{}^\ome\!2$ by setting
$$\buildrul\sim\under\sig
(j)=H_j\left(\buildrul\sim\under\tau\cap
{}^{m_{j\ge}}\!2\right)\ .$$

\subheading{Claim 3.9} $\llvdash_Q$ ``
$\buildrul\sim\under\sig$ is a Cohen real over
$V$".

\pr It is enough to show the following:\hb
$\textstyle\astcirc$\quad for every $\calE>0$  and 
$\eta^*\in {}^{\ome >}\!2$ the following holds:\hb
for every $j<\ome$ large enough and $\nu\in{}^{\ome>}\!
2$  of the length $>j$ there are $q,\langle q_\ell\mid\ell
<\ell^*\rangle$ and $\langle v_\ell\mid\ell
<\ell^*\rangle$  such that
\item{(a)} $q,\langle q_\ell\mid\ell
<\ell^*\rangle$ are in Random 
\smallskip
\item{(b)} $q=\bigcup\limits_{\ell <\ell^*}q_\ell$
\smallskip
\item{(c)} $Lb(q)\ge 1-\calE$
\smallskip
\item{(d)} $v_\ell\trianglelefteq \eta^*$ for
every $\ell <\ell^*$
\smallskip
\item{(e)} $(q_\ell,v_\ell)\llvdash_Q\tagg
\buildrul\sim\under\sig\rhookup [j,\ {\rm length}\ \nu)=
\nu\rhookup [j,\ {\rm length}\ \nu)\tagg$.

\subheading{Proof of $\fakebold{$\astcirc$}$}
Pick $j<\ome$ such that $n_j >\ {\rm length}\
\eta^*$  and $2^{-j}<{\calE\over 2}$. Let $\nu$
be given.  We choose by induction on $k\in [j,\
{\rm length}\ \nu)$ a set $a_k$  and $\langle
q_\eta\mid\eta\in a_k\rangle$  such that 
\item{(a)} $a_k\subseteq{}^{n_k}2$  is nonempty. 
\smallskip
\item{(b)} $a_j$ is a singleton extending
$\eta^*$  
\smallskip
\item{(c)} $\forall\eta\in a_{k+1} (\eta\rhookup
n_k\in a_k)$  and $\forall \eta \in a_k$
$\exists \eta'\in a_{k+1}$ $(\eta\triangleleft
\eta')$.
\smallskip
\item{(d)} for every $\eta\in a_k$
$$(q_\eta,\eta)\llvdash\tagg\bigwedge_{\ell
=j}^{k-1}H_\ell\left(\buildrul\sim\under\tau
\cap{}^{m_\ell\ge}\!2\right)=\nu(\ell)\tagg$$
i.e. $(q_\eta,\eta)\llvdash\tagg\buildrul\sim
\under\sig\rhookup[j,k-1)=\nu\rhookup[ j,k-1)\tagg$
\smallskip
\item{(e)} for $\eta\in a_j$ $q_\eta ={}^\ome\!2$
\smallskip
\item{(f)} for $\eta\in a_k$  $\langle
q_\rho\mid\eta\triangleleft\rho\in a_{k+1}\rangle$
is an antichain of Random above $q_\eta$ and
$\sum\{Lb(q_\rho)\mid\eta\triangleleft\rho\in
a_{k+1}\}\big/Lb(q_\eta)\ge 1-{1\over 2^k}$.

There is no problem in caring on this induction.
This completes the proof of $\astcirc$ and hence
also the theorem.\hfill$\bigsquare$

\vfill\eject
\references {68}
\smallskip
\ref{[Fr]} D. Fremlin, Real-valued measurable
cardinals in Set Theory of the Reals H.~Judah ed.,
Israel Math. Conf. Proceedings (1993), 151-305.

\smallskip
\ref{[Fr2]} D. Fremlin, Measure Algebras, in
J.D. Monk ed., Handbook of Boolean Algebras, North-Holland
(1989), 876-980.

\smallskip
\ref{[Gi-Sh1]} M. Gitik and S. Shelah, On simple
forcing notions and forcing with ideals, Israel
J. Math. 68(2) (1989) 129-160.

\smallskip
\ref{[Gi-Sh2]} M. Gitik and S. Shelah, More on
simple forcing notions and forcing with ideals
Ann. of Pure and Appl. Logic 59 (1993),
219-238. 

\smallskip
\ref{[Pr]} K. Prikry, Ideals and powers of
cardinals, Bull. AMS 81 (1975), 907-909.

\smallskip
\ref{[So]} R. Solovay, Real-valued measurable
cardinals, in: D. Scott, ed. Axiomatic Set
Theory, Proc. Symp. Pure Math. 13(1), (1970), 397-428.  

\smallskip
\ref{[Sh1]} S. Shelah, Cardinal and Arithmetic,
Oxford Logic Guides 29, Oxford Science
Publications 1994.

\smallskip
\ref{[Sh430]} S. Shelah, Further Cardinal Arithmetic,
[Sh430] 

\smallskip
\ref{[Sh460]} S. Shelah, The Generalized Continuum
Hypothesis revisited, [Sh460], to appear.

\smallskip
\ref{[Sh580]} S. Shelah, On Strong Covering, [Sh580].

\end